\documentclass{birkjour}
\usepackage{amsmath, amsthm, amssymb, amsopn, amscd, mathtools, array, enumerate, amsfonts, graphics, hyperref, wasysym, chngcntr, linguex, graphicx, ytableau,youngtab,young,tikz,siunitx}
\usepackage[table]{xcolor}
\usepackage[ruled,vlined]{algorithm2e}
\usepackage{orcidlink}
\newtheorem{theorem}{Theorem}[section]
\newtheorem{lemma}[theorem]{Lemma}
\newtheorem{proposition}[theorem]{Proposition}

\newtheorem{exam}[theorem]{Example}

\newcommand{\field}[1]{\mathbb{#1}}
\newcommand{\Z }{\field{Z}}
\newcommand{\N }{\field{N}}

\DeclareMathOperator{\G}{G}
\DeclareMathOperator{\g}{g}
\DeclareMathOperator{\F}{F}
\DeclareMathOperator{\Co}{C}
\DeclareMathOperator{\pf}{PF}

\DeclareMathOperator{\maxx}{Maximals}
\DeclareMathOperator{\ty}{t}

\begin{document}

%-------------------------------------------------------------------------
% editorial commands: to be inserted by the editorial office
%
%\firstpage{1} \volume{228} \Copyrightyear{2004} \DOI{003-0001}
%
%
%\seriesextra{Just an add-on}
%\seriesextraline{This is the Concrete Title of this Book\br H.E. R and S.T.C. W, Eds.}
%
% for journals:
%
%\firstpage{1}
%\issuenumber{1}
%\Volumeandyear{1 (2004)}
%\Copyrightyear{2004}
%\DOI{003-xxxx-y}
%\Signet
%\commby{inhouse}
%\submitted{March 14, 2003}
%\received{March 16, 2000}
%\revised{June 1, 2000}
%\accepted{July 22, 2000}
%
%
%
%---------------------------------------------------------------------------
%Insert here the title, affiliations and abstract:
%

\title[A Symmetrization for Almost Symmetric Numerical Semigroups]{A Symmetrization for\\ Almost Symmetric Numerical Semigroups}

%----------Author 1
\author[Manuel Delgado]{Manuel Delgado \orcidlink{0000-0002-1935-9908}}
\address{%
CMUP--Centro de Matemática da Universidade do Porto\\
Departamento de Matemática\\
Faculdade de Ciências\\
Universidade do Porto\\
Rua do Campo Alegre s/n\\
4169–007 Porto\\
Portugal}
\email{mdelgado@fc.up.pt}
\thanks{The first author was partially supported by CMUP, member of LASI, which is financed by national funds through FCT – Fundação para a Ciência e a Tecnologia, I.P., under the projects with reference UID/00144/2025 and associated DOI given by https://doi.org/10.54499/UID/00144/2025.}

%----------Author 2
\author[Meral Süer]{Meral Süer \orcidlink{0000-0002-5512-4305}}
\address{Department of Mathematics\\
Faculty of Science and Letters\\
Batman University\\
Batman\\
Türkiye}
\email{meral.suer@batman.edu.tr}

%----------Author 3
\author[Mehmet Yeşil]{Mehmet Yeşil \orcidlink{0000-0002-5220-5177}}
\address{Department of Mathematics\\
	Faculty of Science and Letters\\
	Batman University\\
	Batman\\
	Türkiye}
\email[corresponding author]{mehmet.yesil@batman.edu.tr}

%----------classification, keywords, date


\date{October 29, 2026}
%----------additions
%\dedicatory{}
%%% ----------------------------------------------------------------------

\begin{abstract}
We characterize almost symmetric numerical semigroups in terms of their associated Young diagrams, showing that almost symmetry is equivalent to a simple numerical identity relating the genus, the type, and the number of small elements, and also to a combinatorial condition on the hook lengths of the diagram. We then introduce a symmetrization procedure that associates with each almost symmetric numerical semigroup a symmetric numerical set. This yields a map from the set of almost symmetric numerical semigroups to the set of symmetric numerical sets. We also present experimental data on the image of this map and on the proportion of its elements that are numerical semigroups. The results show, in particular, that the map is far from injective and that its image contains many numerical sets that are not semigroups.
\end{abstract}

%%% ----------------------------------------------------------------------
\maketitle
%%% ----------------------------------------------------------------------
%\tableofcontents
\section{Introduction}

Throughout $\Z$ will denote the set of integers, $\N_{0}$ will denote the set of non-negative integers and $\N$ will denote the set of positive integers. A subset $S$ of $\N_{0}$ that contains zero and has finite complement in $\N_{0}$ is called a \emph{numerical set}. A \emph{numerical semigroup} is a numerical set which is also closed under addition. Note that $\N_{0}$ is a numerical set having empty complement in $\N_0$, and a numerical set is proper if it is not $\N_0$.

Numerical semigroups form a classical and actively studied class of objects at the intersection of commutative algebra, combinatorics, and algebraic geometry. A standard reference is the monograph by Rosales and García-Sánchez \cite{RG}. Symmetry is one of the most fundamental notions in their theory. Symmetric numerical semigroups are precisely those numerical semigroups whose associated semigroup ring is Gorenstein, and they have been studied extensively from both the algebraic and the combinatorial points of view.

Almost symmetric numerical semigroups, introduced by Barucci and Fröberg \cite{BF} in the context of one-dimensional almost Gorenstein rings, generalize the symmetric ones: every symmetric numerical semigroup is almost symmetric of type $1$, while almost symmetric semigroups of higher type form a strictly larger class. A well-known characterization, due to Nari \cite{nari}, relates almost symmetry to the behaviour of the pseudo-Frobenius numbers. Almost symmetric numerical semigroups have since been studied extensively, and they remain a rich source of interesting questions.

A powerful combinatorial tool for studying numerical sets is provided by Young diagrams. A \emph{Young diagram} is a finite collection of boxes arranged in left-justified rows, with the rows ordered from bottom to top, such that the number of boxes in each row is non-decreasing as one moves upwards. Following \cite{AB, KN}, one associates with a proper numerical set a Young diagram by means of a lattice path that encodes the elements and the gaps of the set. This construction establishes a bijection between proper numerical sets and Young diagrams, and it translates arithmetic properties of a numerical set into combinatorial properties of its diagram. This dictionary has been exploited successfully to characterize several classes of numerical semigroups and sets, including Arf semigroups \cite{N1}, symmetric and pseudo-symmetric numerical semigroups \cite{SY} and almost symmetric numerical semigroups \cite{Y}.

The purpose of the present paper is twofold. First, we characterize almost symmetric numerical semigroups in terms of their Young diagrams (Proposition \ref{prop}). Second, we introduce and study a symmetrization procedure for almost symmetric numerical semigroups (Theorem \ref{thm}).

The paper is organized as follows. In Section \ref{sec2}, we recall the basic definitions of numerical sets, numerical semigroups, and Young diagrams, together with the correspondence between numerical sets and Young diagrams. In Section \ref{sec3}, we characterize almost symmetric numerical semigroups via their Young diagrams. In Section \ref{sec4}, we introduce the symmetrization construction and prove that it produces symmetric numerical sets. In Section \ref{sec5}, we describe an algorithm for finding the symmetrization of a given almost symmetric numerical semigroup and present experimental results, together with a discussion of the injectivity and surjectivity of the symmetrization function.

\section{Young Diagrams and Numerical Semigroups}\label{sec2}

Consider a Young diagram $Y$ with $n$ columns and $g$ rows. The number of boxes in a column (or row) is called the length of that column (or row). The hook of a box in $Y$ is the shape formed by the boxes directly to the right of it, the boxes directly below it, and the box itself. The number of boxes in the hook of a box is called the hook length of that box. The standard reference for Young diagrams is the monograph of Fulton \cite{F}.

\begin{exam}\label{ex1}
	The figure below shows a Young diagram with $7$ columns and $9$ rows. The boxes containing bullets form the hook of the box in the second column (from left to right) and second row (from top to bottom).
	\begin{center}
		$\tiny{\ytableaushort{\none,\none\bullet\bullet\bullet, \none\bullet,\none\bullet,\none\bullet,\none\bullet}*{7,4,4,2,2,2,1,1,1}}$
	\end{center}
	The hook length of that box is $7$.
\end{exam}

Henceforth, we assume that $S$ is a proper numerical set. The gaps of $S$ are positive integers not in $S$, and $\G(S)$ stands for the set of those. The genus of $S$, denoted by $\g(S)$, is the cardinality of $\G(S)$. The Frobenius number of $S$, denoted by $\F(S)$, is the largest gap of $S$. Furthermore, $\Co(S):=\F(S)+1$ is called the conductor of $S$. Equivalently, $\Co(S)$ is the smallest element of $S$ such that every integer greater than or equal to $\Co(S)$ belongs to $S$. Conventionally, we note that $\F(\N_{0})=-1$ and $\Co(\N_{0})=0$.

The elements of $S$ that are less than or equal to the conductor are said to be the small elements. If $S$ contains $n$ non-zero small elements arranged in ascending order and listed as $0 = s_0 < s_1 < \dots < s_{n}$, we represent $S$ as follows: $$S=\{0,s_{1},\dots,s_{n}=\Co(S),\rightarrow\}.$$ Here, the arrow at the end means that all integers greater than $\Co(S)$ are in $S$.

For a non-negative integer $x$, we define the set difference $S-x$ as follows: $$S-x = \{s-x \mid s \in S \text{ and } s \geq x\}\cup\{0\}.$$ If $x\geq \F(S)$, then $S-x=\N_0$. If $x<\F(S)$, then $S-x$ forms a numerical set with gaps as follows; $$\G(S-x) = \{g-x \mid g \in \G(S) \text{ and } g > x\}.$$ Consequently, the Frobenius number of $S-x$ is given by $\F(S-x) = \F(S) - x$. However, if $x$ is equal to a small element $s_i$ within $S$, then $S-x$ can be expressed as: $$S-x=\{ 0=s_i-x,s_{i+1}-x,\dots,s_n-x,\rightarrow\}.$$ 

A gap $z\in\G(S)$ is called a pseudo-Frobenius number of $S$ if $z+s\in S$ for every non-zero element $s\in S$. We denote the set of pseudo-Frobenius numbers of $S$ by $\pf(S)$, and the cardinality of this set by $\ty(S)$, which is called the type of $S$. For $z,w\in\Z$, we define a partial order on $\Z$ by setting $z\leq_{S}w$ if $w-z\in S$. For a subset $A$ of $S$, the set of maximal elements of $A$ with respect to the relation $\leq_{S}$ is denoted by $\maxx_{\leq_S}(A)$. By the definition, we observe that pseudo-Frobenius numbers are the maximal elements in $\G(S)$ with respect to $\leq_S$, that is $\pf(S)=\maxx_{\leq_S}(\G(S))$. Then we can easily deduce that $z\in\G(S)$ if and only if $f-z\in S$ for some $f\in\pf(S)$.

A numerical set $S$ is called symmetric if for every $a\in[0,\F(S)]$ exactly one of $a$ or $\F(S)-a$ is an element of $S$. When $S$ is a semigroup, being symmetric is equivalent to the property that $2\g(S)=\F(S)+1$.

More generally, it is well-known that $2\g(S)\geq \F(S)+\ty(S)$ is valid for any numerical semigroup $S$. When the equality holds, i.e. $2\g(S)= \F(S)+\ty(S)$, $S$ is called almost symmetric. This means that every symmetric numerical semigroup is almost symmetric of type $1$.

For a numerical set $S$, we can create a Young diagram, labelled as $Y_S$, that corresponds to $S$ by drawing a continuous polygonal path starting at the origin in $\Z^2$. We begin with $s=0$, and for each subsequent $s$:
\begin{enumerate}
	\item if $s \in S$, we draw a unit-length line to the right.
	\item if $s \notin S$, we draw a unit-length line upwards.
\end{enumerate}
We continue this process until we reach $s=\F(S)$. The region bounded by this polygonal path and the horizontal line at height $\g(S)$ defines the corresponding Young diagram $Y_S$. It is clear that every Young diagram corresponds to a unique proper numerical set. Consequently, the mapping from $S$ to $Y_S$ establishes a one-to-one correspondence, i.e. a bijection, between the set of proper numerical sets and the set of Young diagrams.

When we have a numerical set $S=\{0,s_{1},\dots,s_{n-1},s_{n}=\Co(S),\rightarrow\}$ with the corresponding Young diagram $Y_{S}$, by the construction process, it is evident that $Y_{S}$ has $\g(S)$ rows and $n$ columns. For each $i=0,1,\dots,n-1$, we denote by $\G_i(S)$ the set of hook lengths of the boxes in the $i$th column (from left) of \(Y_S\), which corresponds to $s_i$. Also, the $i$th row (from bottom) corresponds to the $i$th gap of $S$. The hook length of the box in the first column and $i$th row is the $i$th gap of $S$. Thus, $\G_{0}(S)=\G(S)$. Furthermore, we denote the set of hook lengths of boxes in the top row of $Y_S$ by $\G'(S)$.

\begin{proposition} \label{prp} \cite[Section 2]{N1}
	Let $S=\{0,s_{1},\dots,s_{n-1},s_{n},\rightarrow\}$ be a numerical set with corresponding Young diagram $Y_{S}$. Then:
	\begin{enumerate}
		\item $\G'(S)=\{\F(S)-s_{0},\F(S)-s_{1},\dots,\F(S)-s_{n-1}\}$,
		\item $S$ is a numerical semigroup if and only if $\G_{i}(S)\subseteq \G(S)$ for each $i\in \{ 0,1,\dots,n-1 \}$.
	\end{enumerate}
\end{proposition}

\begin{lemma} \label{lem}
	Let $S=\{0,s_{1},\dots,s_{n-1},s_{n},\rightarrow\}$ be a numerical semigroup with corresponding Young diagram $Y_S$. Then $\pf(S)\setminus \{\F(S)\}\subseteq \G(S)\setminus \G'(S)$.
\end{lemma}

\begin{proof}
	By Proposition \ref*{prp} (1), we know that $\G'(S)=\{\F(S)-s_{0},\F(S)-s_{1},\dots,\F(S)-s_{n-1}\}$, and by Proposition \ref*{prp} (2), it follows immediately that $\G'(S)\subseteq\G(S)$. Since we have $ s_i=\F(S)-(\F(S)-s_i)\in S$ for each $i\in \{ 0,1,\dots,n-1 \}$, $\maxx_{\leq_S}(\G'(S))=\{ \F(S)\}$. Hence, the only pseudo-Frobenius number in $\G'(S)$ is $\F(S)$.
\end{proof}

\section{Young Diagrams of Almost Symmetric Numerical Semigroups} \label{sec3}
In this section, we characterize almost symmetric numerical semigroups in terms of their Young diagrams.

\begin{theorem} \cite[Section 2]{nari}
	Let $S=\{0,s_{1},\dots,s_{n-1},s_{n},\rightarrow\}$ be a numerical semigroup with $\pf(S)=\{ f_1,\dots,f_t=\F(S)\}$. Then the following are equivalent.
	\begin{enumerate}
		\item $S$ is almost symmetric,
		\item $a\in \Z\setminus S$ implies that either $a\in\pf(S)$ or $\F(S)-a\in S$,
		\item $f_i+f_{t-i}=\F(S)$ for all $i\in\{ 1,\dots,t-1 \}$.
	\end{enumerate}
	
\end{theorem}

\begin{proposition} \label{prop}
	Let $S=\{0,s_{1},\dots,s_{n-1},s_{n},\rightarrow\}$ be a numerical semigroup with corresponding Young diagram $Y_S$. Then the following are equivalent.
	\begin{enumerate}
		\item $S$ is almost symmetric,
		\item $\g(S)-\ty(S)=n-1$,
		\item $\G(S)\setminus \G'(S)=\pf(S)\setminus\{\F(S)\}$.
	\end{enumerate}
\end{proposition}

\begin{proof}
	$(1)\implies(2)$ Suppose that $S$ is almost symmetric, i.e. $2\g(S)= \F(S)+\ty(S)$. Since the conductor of $S$ is the sum of genus and number of small elements, $\F(S)=\g(S)+n-1$. Then we get $2\g(S)=\g(S)+n-1+\ty(S) \implies \g(S)-\ty(S)=n-1.$
	
	$(2)\implies(3)$ Assume that $\g(S)-\ty(S)=n-1$. This means that $\g(S)-n=\ty(S)-1$. Since $|\G(S)|=\g(S)$, $|\G'(S)|=n$ and $\G'(S)\subseteq\G(S)$, we get $|\G(S)\setminus\G'(S)|=\g(S)-n$. On the other hand $|\pf(S)\setminus\{\F(S)\}|=\ty(S)-1$. Now by Lemma \ref{lem}, it follows that  $\g(S)-n=\ty(S)-1\implies\G(S)\setminus \G'(S)=\pf(S)\setminus\{\F(S)\}$.
	
	$(3)\implies(1)$ Assume now that $\G(S)\setminus \G'(S)=\pf(S)\setminus\{\F(S)\}$. It follows that $\g(S)-n=\ty(S)-1$, and so $\g(S)=n-1+\ty(S)$. Then adding $\g(S)$ to both sides of the latter equation gives $2\g(S)=\g(S)+n-1+\ty(S)=\F(S)+\ty(S)$. Hence, $S$ is almost symmetric.
\end{proof}

This result, stated above, provides a purely combinatorial criterion for almost symmetry and complements the existing characterizations in the literature.

\begin{exam} \label{ex2}
	The numerical set $S=\{0, 4, 8, 9, 12, 13, 14, 16,\rightarrow \}$ is a numerical semigroup with $\G(S) = \{1,2,3,5,6,7,10,11,15\}$. Here $S$ has $7$ non-zero small elements and $9$ gaps. Therefore, the corresponding Young diagram $Y_{S}$ has $7$ columns and $9$ rows, and the hook lengths of boxes in the first column and the first row of $Y_{S}$ are shown in the picture below.
	\begin{center}
		$\tiny{\ytableaushort{{15}{11}76321,{11},{10},7,6,5,3,2,1}*{7,4,4,2,2,2,1,1,1}}$
	\end{center}
	We can easily find $\pf(S)=\{ 5,10,15 \}$ and $\ty(S)=3$. Since $2\g(S)=\F(S)+\ty(S)$, $S$ is actually almost symmetric. Using its corresponding Young diagram, we can also see that $\G(S)\setminus\G'(S)=\{5,10\}=\pf(S)\setminus\{\F(S)\}$.
\end{exam}

\section{A Symmetrization for Almost Symmetric Numerical Semigroups via Young Diagrams} \label{sec4}

In this section, we construct symmetric numerical sets from almost symmetric numerical semigroups using their Young diagrams.

Let $S=\{0,s_{1},\dots,s_{n-1},s_{n},\rightarrow\}$ be a numerical semigroup with $\G(S)=\{ a_1,\dots,a_g=\F(S) \}$ and the corresponding Young diagram $Y_S$. As we know the $i$th row (from bottom) of $Y_S$ corresponds to the $i$th gap of $S$, we denote that row by $R_{a_i}$. If we delete $R_{a_i}$ from $Y_S$ for some $i\in\{1,\dots,g\}$, we get a new Young diagram $Y_S'$ and also its corresponding numerical set $S'$. When we delete $R_{a_i}$ from $Y_S$, the number of rows decreases by one, while the number of columns remains unchanged, which means that $S'$ has $n$ small elements and $g-1$ gaps. Moreover, deleting $R_{a_i}$ from $Y_S$ affects only the gaps and small elements of $S$ greater than $a_i$. So, the resulting numerical set $S'$ can be obtained as follows:
$$S'=\{ 0,s_1',\dots,s_n',\rightarrow \}
\text{ where } s_j'=\begin{cases}
	s_j & \text{ if } s_j<a_i \\
	s_j-1 & \text{ if } s_j>a_i
\end{cases}$$ $$ \text{ with } \G(S')=\{ a_1',\dots,a_{g-1}' \} \text{ where } a_j'=\begin{cases}
	a_j & \text{ if } a_j<a_i \\
	a_{j+1}-1 & \text{ if } a_j\geq a_i
\end{cases}.$$
Hence, if $s_j$ is the smallest element of $S$ that is greater than $a_i$, then we get $$S'=\{ 0,s_1,\dots,s_{j-1},s_j-1,\dots,s_n-1,\rightarrow \}$$
and $$\G(S')=\{ a_1,\dots,a_{i-1},a_{i+1}-1,\dots,a_g-1\}.$$
We can also consider deleting more than one row from $Y_S$. In this case, $S'$ can be obtained by applying the preceding procedure inductively.

\begin{exam}
	Let $S=\{0, 4, 8, 9, 12, 13, 14, 16,\rightarrow \}$ with $\G(S)=\{ 1,2,3,5,6,7,10,11,15 \}$. If we delete $R_{5}$ and $R_{10}$ from $Y_S$, we get $S'=\{ 0,4,7,8,10,11,12,14,\rightarrow\}$ and $\G(S')=\{ 1,2,3,5,6,9,13 \}$.
	\tiny{\begin{center}
			\begin{tabular}{ccccc}
				\ytableaushort{{15},{11},{10},7,6,5,3,2,1}*{7,4,4,2,2,2,1,1,1}*[*(red)]{0,0,4,0,0,2,0,0,0}&$\longrightarrow$&\ytableaushort{{13},9,6,5,3,2,1}*{7,4,2,2,1,1,1}\\
				$\underbrace{\hspace{11.2em}}_{\displaystyle Y_S}$&&$\underbrace{\hspace{11.2em}}_{\displaystyle Y_{S}'}$
			\end{tabular}
	\end{center}}
	
\end{exam}

In particular, we focus on deleting rows of $Y_S$ which correspond to pseudo-Frobenius numbers of $S$ that are smaller than the Frobenius number. The following theorem gives this characterization for almost symmetric numerical semigroups which extends \cite[Lemma 14]{Y} to the general case.

\begin{theorem} \label{thm}
	Let $S=\{0,s_{1},\dots,s_{n-1},s_{n},\rightarrow\}$ be an almost symmetric numerical semigroup of type $t$, with pseudo-Frobenius numbers $f_1<\cdots<f_t=\F(S)$. Let $Y_S$ be the Young diagram corresponding to $S$ and $Y_S'$ be the Young diagram obtained by deleting each $R_{f_i}$ for $i\in\{1,\dots,t-1\}$ from $Y_S$. Then the numerical set $S'$ corresponding to $Y_S'$ is symmetric.
\end{theorem}

\begin{proof}
	By applying the row-deletion procedure inductively, we obtain $S'=\{ 0,s_1',\dots,s_n',\rightarrow \}$ where
	$$ s_j'=\begin{cases}
		s_j & \text{ if } s_j<f_1 \\
		s_j-i & \text{ if } f_i<s_j<f_{i+1}\\
		s_j-t+1 & \text{ if } f_{t-1}<s_j
	\end{cases}$$
	and $\G(S')=\{ a_1',\dots,a_{g-t+1}'\}$  where $a_{j}'=a_j$ if $a_j<f_1$ and 
	$$a_{j-i}'=\begin{cases}
		a_{j}-i & \text{ if } f_i<a_j<f_{i+1}\\
		a_{j}-t+1 & \text{ if } f_{t-1}=a_j
	\end{cases}$$
	for each $i=1,\dots,t-1$, and we get $\F(S')=\F(S)-t+1$. Take $a\in[0,\F(S')]\cap\N$. Assume that $a\notin S'$. If $a=\F(S')$, $\F(S')-a=0\in S'$. We may therefore assume that $a\neq\F(S')$. Hence, $a=a_j'\in\G(S')$ for some $j=1,\dots,g-t$. Then there exist $a_j\in\G(S)$ such that $a_j'=a_j-i$ and $f_i<a_j<f_{i+1}$ for some $i=1,\dots,t-1$. Since $a_j\in\G(S)\setminus\pf(S)$ and $S$ is almost symmetric, $\F(S)-a_j\in S$. Furthermore, $\F(S)-f_{i+1}<\F(S)-a_j<\F(S)-f_i$. This means we have $f_{t-i-1}<\F(S)-a_j<f_{t-i}$. Thus, we get $$\F(S)-a_j-t+i+1\in S' \implies \F(S')-(a_j-i)=\F(S')-a_j'=\F(S')-a\in S'.$$
	Hence, $S'$ is symmetric.
\end{proof}

\section{Experimental Results}\label{sec5}

We begin this section by translating the construction in Theorem \ref{thm} into a pseudocode. This translation appears in Algorithm \ref{alg:Symmetrization} and can be converted into a GAP \cite{GAP} function immediately. We then present some data to compare the cardinalities of the sets of almost symmetric numerical semigroups and symmetric numerical sets.

\begin{algorithm}[h]  \caption{Symmetrization}\label{alg:Symmetrization}
	% \DontPrintSemicolon
	\SetKwComment{cgap}{\# }{} %comment in gap style
	\KwData{An almost symmetric numerical semigroup $S$}
	\KwResult{The symmetric numerical set $S'$}
	
	$se := \mathrm{SmallElements}(S)$\\
	$pf:= \mathrm{PseudoFrobenius}(S)$\\
	$tp:= Type(S)$\\
	$S':= \mathrm{Filtered}(se, s \,-\!\!> s < pf[1])$\\
	\For{$i$  in $[1..tp-1]$}{
		$S':= \mathrm{Concatenation}(S',\mathrm{Filtered}(se, s \,-\!\!> s < pf[i+1]))$  and  $s > pf[i])-i)$}
	$\mathrm{Add}(S',se[\mathrm{Length}(se)]-tp+1)$\\
	\KwRet $S'$
\end{algorithm}

Assuming that Algorithm \ref{alg:Symmetrization} is implemented in GAP under the name \texttt{Symmetrization} and that the \texttt{NumericalSgps} package \cite{DS} is loaded, the data in Table \ref{tab} can be computed by executing the following code with the appropriate value $F$. (Warning: the current implementations require a considerable amount of memory, and may not work for large $F$.)\\

\noindent
\texttt{asns := AlmostSymmetricNumericalSemigroupsWithFrobeniusNumber(F);}\\
\texttt{setns := Set(asns, s -> Symmetrization(s));}\\
\texttt{l := List([1..Length(setns)], i ->} \\
\hspace*{2em}\texttt{RepresentsSmallElementsOfNumericalSemigroup(setns[i]));}\\
\texttt{[F, Length(asns),Length(setns), SizeBlist(l)]};\\

Given an almost symmetric numerical semigroup, the Algorithm \ref{alg:Symmetrization} returns a symmetric numerical set, and the construction is well-defined. Consequently, we obtain a map from the set of almost symmetric numerical semigroups to the set of symmetric numerical sets. We call this map the \emph{symmetrization function} and denote it by $\Phi$.

Let $\mathcal{AS}$ be the set of almost symmetric numerical semigroups and $\mathcal{S}$ be the set of symmetric numerical sets. The symmetrization function $\Phi:\mathcal{AS}\longrightarrow\mathcal{S}$ is not injective. For instance, if
\[
S=\{0,4,6,8,\rightarrow\},\quad T=\{0,5,6,8,\rightarrow\},
\]
then $S,T\in\mathcal{AS}$ and $\Phi(S)=\Phi(T)=\{0,3,4,6,\rightarrow\}$.

For a positive integer $F$, let $\mathcal{AS}_F$ denote the set of almost symmetric numerical semigroups with Frobenius number $F$. Let $\mathcal{S}_F$ denote the image of $\mathcal{AS}_F$ under $\Phi$, i.e. the set of symmetric numerical sets obtained by applying $\Phi$ to elements of $\mathcal{AS}_F$. Since the output of $\Phi$ is a numerical set but not necessarily a semigroup, a natural question is: how many elements of $\mathcal{S}_F$ are actually numerical semigroups? To this end, we denote by $\mathcal{NS}_F$ the subset of $\mathcal{S}_F$ consisting of those symmetric numerical sets that are also semigroups.

The experimental results in Table \ref{tab} compare the cardinalities of $\mathcal{AS}_F$, $\mathcal{S}_F$, and $\mathcal{NS}_F$. This ratio grows rapidly, indicating that the symmetrization map collapses a large number of almost symmetric numerical semigroups onto a comparatively small number of symmetric numerical sets.

\begin{table}[h]
	\centering
	\setlength{\tabcolsep}{10pt}
	\setlength{\arrayrulewidth}{0.5mm}
	{\rowcolors{2}{black!10}{}
		\begin{tabular}{| r  r  r  r  |  r  r  r  r  |}
			\hline
			\(F\)&\( |\mathcal{AS}_F| \)&\( |\mathcal{S}_F|  \) &\( |\mathcal{NS}_F|  \) &
			\(F\)&\( |\mathcal{AS}_F| \)&\( |\mathcal{S}_F|  \) &\( |\mathcal{NS}_F|  \)                                           
			\\ \hline		
			1	&1	&1	&1	&40	&\num{7694}	&558	&372\\
			2	&1	&1	&1	&41	&\num{14966}	&\num{1228}	&862\\
			
			3	&2	&2	&2	&42	&\num{10341}	&666	&466\\
			4	&2	&2	&2	&43	&\num{22600}	&\num{1659}	&\num{1179}\\
			
			5	&4	&4	&4	&44	&\num{17607}	&\num{1016}	&736\\
			6	&3	&3	&3	&45	&\num{31918}	&\num{1832}	&\num{1314}\\
			
			7	&7	&6	&6	&46	&\num{24768}	&\num{1467}	&\num{853}\\
			8	&6	&5	&5	&47	&\num{51435}	&\num{2882}	&\num{1945}\\
			
			9	&11	&8	&8	&48	&\num{38043}	&\num{1522}	&\num{994}\\
			10	&9	&7	&6	&49	&\num{77177}	&\num{3760}	&\num{2486}\\
			
			11	&20	&13	&12	&50	&\num{56054}	&\num{2432}	&\num{1542}\\
			12	&15	&8	&8	&51	&\num{108614}	&\num{4459}	&\num{2948}\\
			
			13	&32	&18	&17	&52	&\num{89724}	&\num{3151}	&\num{1898}\\
			14	&24	&13	&13	&53	&\num{174024}	&\num{6820}	&\num{4329}\\
			
			15	&47	&21	&20	&54	&\num{120550}	&\num{3778}	&\num{2304}\\
			16	&42	&18	&15	&55	&\num{259874}	&\num{8665}	&\num{5538}\\
			
			17	&83	&34	&30	&56	&\num{201926}	&\num{5593}	&\num{3359}\\
			18	&58	&21	&19	&57	&\num{365294}	&\num{10219}	&\num{6244}\\
			
			19	&131	&48	&43	&58	&\num{282173}	&\num{7917}	&\num{4115}\\
			20	&103	&30	&30	&59	&\num{581512}	&\num{15704}	&\num{9469}\\
			
			21	&188	&55	&49	&60	&\num{427069}	&\num{7708}	&\num{4729}\\
			22	&151	&44	&35	&61	&\num{866400}	&\num{20734}	&\num{12185}\\
			
			23	&323	&84	&73	&62	&\num{628533}	&\num{13801}	&\num{7541}\\
			24	&241	&48	&42	&63	&\num{1219787}	&\num{23439}	&\num{14005}\\
			
			25	&494	&109	&93	&64	&\num{997083}	&\num{17277}	&\num{8879}\\
			26	&364	&79	&66	&65	&\num{1915186}	&\num{34142}	&\num{19497}\\
			
			27	&718	&138	&116	&66	&\num{1329239}	&\num{19947}	&\num{10490}\\
			28	&605	&99	&75	&67	&\num{2855488}	&\num{47581}	&\num{26500}\\
			
			29	&\num{1192}	&213	&174	&68	&\num{2199659}	&\num{29242}	&\num{15878}\\
			30	&821	&111	&91	&69	&\num{4001514}	&\num{54115}	&\num{30054}\\
			
			31	&\num{1827}	&286	&229	&70	&\num{3071586}	&\num{39189}	&\num{18695}\\
			32	&\num{1441}	&181	&147	&71	&\num{6297210}	&\num{80923}	&\num{43416}\\
			
			33	&\num{2599}	&329	&261	&72	&\num{4629706}	&\num{43416}	&\num{22051}\\
			34	&\num{2034}	&257	&177	&73	&\num{9327934}	&\num{106317}	&\num{55771}\\
			
			35	&\num{4243}	&477	&373	&74	&\num{6755646}	&\num{71591}	&\num{34458}\\
			36	&\num{3172}	&272	&207	&75	&\num{13040617}	&\num{118432}	&\num{63341}\\
			
			37	&\num{6496}	&690	&510	&76	&\num{10603230}	&\num{86922}	&\num{40702}\\
			38	&\num{4725}	&457	&329	&77	&\num{20481464}	&\num{182931}	&\num{92213}\\
			
			39	&\num{9166}	&780	&597	&78	&\num{14123754}	&\num{102019}	&\num{47682}\\
			\hline	                  
	\end{tabular}}
	\caption{\label{tab}}
\end{table}

\begin{figure}[h]
	\centering
	\includegraphics[width=0.75\textwidth]{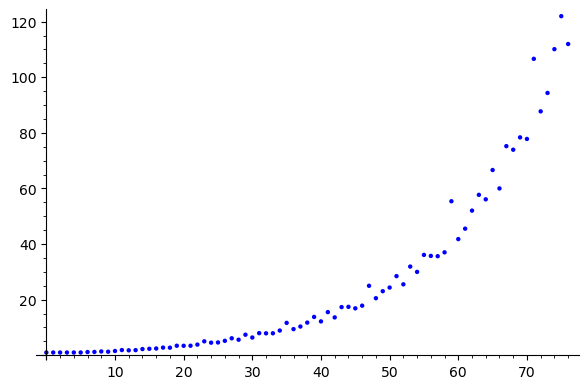}
	\caption{The ratio \(|\mathcal{AS}_F| / |\mathcal{S}_F|\) for \(1\le F\le 78\).}
	\label{fig:AS_by_S}
\end{figure}

Using this data, we can observe two interesting asymptotic behaviours. First, Figure \ref{fig:AS_by_S} shows the ratio $|\mathcal{AS}_F| / |\mathcal{S}_F|$ as a function of $F$. This ratio grows rapidly, indicating that the symmetrization map collapses a large number of almost symmetric numerical semigroups onto a comparatively small number of symmetric numerical sets.

\begin{figure}[h]
	\centering
	\includegraphics[width=0.75\textwidth]{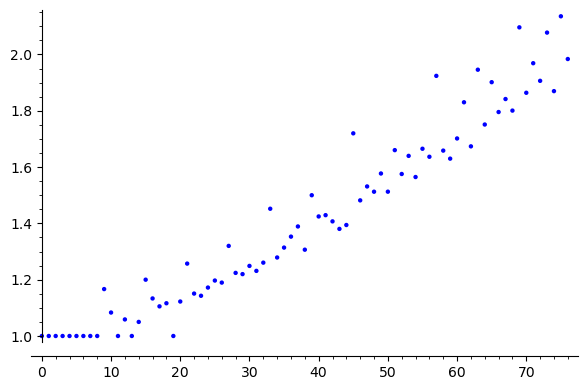}
	\caption{The ratio \(|\mathcal{S}_F| / |\mathcal{NS}_F|\) for \(1\le F\le 78\)}
	\label{fig:S_by_NS}
\end{figure}

Second, Figure \ref{fig:S_by_NS} shows the ratio $|\mathcal{S}_F| / |\mathcal{NS}_F|$. This ratio starts at $1$ for small values of $F$, meaning that, in this range, every symmetric numerical set produced by $\Phi$ is again a numerical semigroup. As $F$ increases, the ratio generally grows, with some oscillation, exceeding $2$ at $F=78$. This suggests that, in the computed range, numerical semigroups account for roughly half of the elements of $\mathcal{S}_F$, while the remaining elements are symmetric numerical sets that are not closed under addition.

Taken together, these two observations suggest that the symmetrization function $\Phi$ is far from being injective, and that its image contains a significant proportion of numerical sets that are not numerical semigroups. A natural question is whether $\Phi$ is surjective, that is, whether every symmetric numerical set arises as the symmetrization of some almost symmetric numerical semigroup. Determining whether $\Phi:\mathcal{AS}\longrightarrow\mathcal{S}$ is onto, or equivalently characterizing its image, would be an interesting direction for future work.

%\subsection*{Acknowledgment}
%Many thanks to our \TeX-pert for developing this class file.

% ------------------------------------------------------------------------
\end{document}